\documentclass[12pt]{iopart}	

\usepackage{iopams,tabularx,amscd,amsthm}
\usepackage[arrow,matrix,curve]{xy}
\usepackage{graphicx}
\usepackage{geometry}

 \newtheorem{thm}{Theorem}

 \theoremstyle{definition}
  \newtheorem{ex}{Example}
 \theoremstyle{remark}

\newcommand{\Sc}[2]{\langle #1,#2\rangle}
\newcommand{\Hom}{\mathrm{Hom}}

\newcommand{\Ext}{\mathrm{Ext}}
\newcommand\Zn{\mathbb{Z}}
\newcommand\Nn{\mathbb{N}}
\newcommand\Cn{\mathbb{C}}

\newcommand{\Sp}{$\textrm{Spec}$}
\newcommand{\Pj}{$\textrm{Proj}$}

\begin{document}

\title{Representations of posets: Linear versus Unitary}

\author{Vyacheslav Futorny$^1$, Yurii Samoilenko$^2$, Kostyantyn Yusenko$^{1,2}$}

\address{$^1$ IME-USP, S\~ao-Paulo, Brasil}

\address{$^2$ Institute of Mathematics, Kiev, Ukraine}

\ead{kay.math@gmail.com}

\begin{abstract}
A number of recent papers treated the representation theory of partially ordered
sets in unitary spaces with the so called orthoscalar relation.
Such theory generalizes the classical theory which studies
the representations of partially ordered sets in linear spaces. It happens that the
results in the unitary case are well-correlated with those in
the linear case. The purpose of this article is to shed light on
this phenomena.

\,


\end{abstract}

\section*{0. Introduction}

In the second part of XX century the new results concerning 
to representation theory in linear spaces were obtained, in particular, 
the results about relative position of several subspaces and 
partially ordered collections of subspaces in linear spaces (see, for example, \cite{GP, N}).
Representations of partially ordered sets (\emph{posets} in
sequel) were first introduced by
L.Nazarova and A.Roiter (see \cite{NazarovaRoiter}), using the language of "matrix problems",
as a tool to prove the second Brauer-Thrall conjecture. Later
P.Gabriel gave an equivalent (up to a finite number of
indecomposables) definition of representations using "subspace"
language (see, for example, \cite[Chapter 3]{Simson1992}). In the
sequel we use the latter terminology.

One can try to develop
a similar theory for unitary spaces. Namely, one can  transfer
straightforward  the definitions  from linear to unitary case
keeping in mind that the morphisms between two representations
should preserve their unitary structure; that is, the morphisms are unitary
maps. With this restriction the classification problem becomes extremely difficult already in
very simple situations. For example the problem of classifying unitary representations of the poset $\mathcal P=\{a_1,a_2,a_3 \ | \ a_2\prec a_3 \}$  is hopeless, since by \cite{KruglyakSamoilenko1,KruglyakSamoilenko2} it contains the problem of classifying any system of operators in unitary spaces (see \cite[Chapter 3]{OstrovskyiSamoilenko} for more details).

The idea, hence, is to consider those representations which
satisfy some additional conditions. In such a way there were
defined orthoscalar representations of  quivers (see
\cite{KrugNazRoi,KrugRoi}) and orthoscalar representations of
posets (see, for example, \cite{AlbeverioOstrSamoilenko,SY}). Certainly the
additional conditions can be chosen differently.
Orthoscalarity relation appears in different areas of mathematics. For example, using \cite{Klyachko,KnuS}, one can show that every indecomposable orthoscalar representation of a poset generates a stable reflexible sheaf on certain toric variety. This is one of the 
motivation to study such representations. 

Analyzing the  results related to orthoscalar representations of
quivers and posets (see
\cite{AlbeverioOstrSamoilenko, KrugNazRoi,KrugRoi,SY}  and
references therein) one can see the mysterious  (on the
first sight) connection with the classical results (more details in Section
1). The main goal of this paper is  to explain this
phenomena.

For a given finite poset $\mathcal P$ we consider the variety $R_{\mathcal
P,d}$ of all representations of $\mathcal P$ with dimension vector
$d=(d_0;d_i)_{i\in \mathcal P}$.
We show that the classes of unitary non-equivalent orthoscalar representations of $\mathcal P$ with
dimension vector $d$ can be viewed as the symplectic quotient of
$R_{\mathcal P,d}$. On the other side we consider the GIT quotient of
$R_{\mathcal P,d}$ with respect to some polarization of this
variety. This quotient can be viewed as the set of isomorphism
classes of polystable linear representations of $\mathcal P$ with
dimension $d$. We show that the identification between these two
quotients is a consequence of the fundamental Kirwan-Ness theorem
(see, for example, \cite{MumfordForganyKirwan}). 

Let us mention that the idea of the symplectic reduction in
similar contexts is not new. For example, it appeared in
\cite{kin} (constructing moduli spaces for quivers) and in \cite{Knu} (in
connection with Horn's problem). In this text we  apply it in
the setting of partially ordered sets.

\subsection{Structure of the paper.}

In Section \ref{s_pre} we recall basic facts about the posets and
their indecomposable representations in linear spaces and indecomposable orthoscalar 
representations in unitary spaces, in particular we
compare some known results about the structure of representations in
both cases. In Section~\ref{s_symplect} we give some basic facts about
symplectic quotients. Section~\ref{s_git} is dedicated
to the GIT quotients. We explain the construction and state the
Kirwan-Ness theorem which connects symplectic and GIT quotients. In
Section~\ref{s_corr} we prove Theorem~\ref{posCor} which relates
representations of posets in linear spaces with the ones in unitary
spaces with orthoscalar condition. In Section~\ref{s_dimension} we
estimate the dimension of the variety of uni-classes of
representations of a given poset with a  given dimension vector
using the connection between posets and bound quivers.


\section{Posets and their
representations.} \label{s_pre}

Let $\mathcal P$ be a finite poset with partial order relation $\prec$. Denote by $\mathcal P^*$ the extension of $\mathcal P$ by a unique maximal element $^*$. To
$\mathcal P$ we associate  the Hasse quiver of $\mathcal P^*$
which is denoted by $Q^{\mathcal P}$. The \textit{order} of $\mathcal P$ is
the number of its elements. A poset is said
to be \emph{primitive} and denoted by $(n_1,\ldots,n_s)$ if it is a disjoint (cardinal) sum of
linearly ordered sets with the orders $n_i$. 
A poset $\mathcal P$ is primitive if and only if the corresponding quiver $Q_\mathcal P$
has a star-shaped form.

\begin{ex}  Let $\mathcal P=(1,2)$ be the poset consisting of three elements
$a_1, a_2, a_3$ with a unique relation $a_2\prec a_3$. Then the
associated quiver $Q^{\mathcal P}$ has the following form:

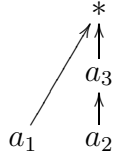
\begin{center}
$\xymatrix @-1pc { &{*}\ar@{<-}[d]\ar@{<-}[ddl]&\\
                   &{a_3}\ar@{<-}[d]&\\
            {a_1}&{a_2}&}$
            \qquad
\end{center}

\end{ex}

\subsection{Linear representations}

Following P.Gabriel (see \cite[Chapter 3]{Simson1992}) we define the category
$\mathcal P-sp$ of representations of $\mathcal P$ in linear spaces.
The objects of $\mathcal P-sp$ are systems $(V;V_i)_{i\in \mathcal
P}$, in which $V$ is a vector space, each $V_i$ is a subspace of $V$ and
$V_i\subset V_j$ if $i\prec j$. The set of morphisms between two
objects $(V;V_i)_{i\in \mathcal P}$ and $(\tilde V;\tilde V_i)_{i\in
\mathcal P}$  consists of the linear maps $f:V\rightarrow \tilde V$ such
that $f(V_i)\subset \tilde V_i$ for all $i\in \mathcal P$. Two
objects are \emph{isomorphic} if there exists an invertible
$f:V\rightarrow \tilde V$ with $f(V_i)=\tilde V_i$ for all $i \in
\mathcal P$. The category $\mathcal P-sp$ is additive with usual
direct sums. An object $(V;V_i)_{i\in \mathcal P}$ is said to be
\emph{indecomposable} it is not isomorphic to a direct sum of two
non-zero objects.

\begin{ex} Let $\mathcal P=(1,1)$ be a poset consisting of two
incomparable elements. Fixing $\lambda \in \mathbb C$ we build the following  linear two-dimensional representation of $\mathcal P$
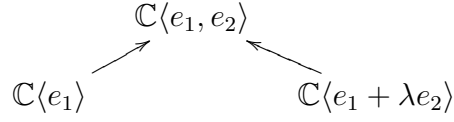
\begin{center}
$\xymatrix @-1pc { &\mathbb C\langle e_1,e_2\rangle\ar@{<-}[dl]\ar@{<-}[dr]&\\
            \mathbb C \langle e_1\rangle&&\mathbb C\langle e_1+\lambda e_2\rangle}$
            \qquad
\end{center}
It is easy to see that any such representation with
$\lambda \neq 0$ is equivalent to the one with $\lambda=1$  and
hence it splits into one-dimensional representations:
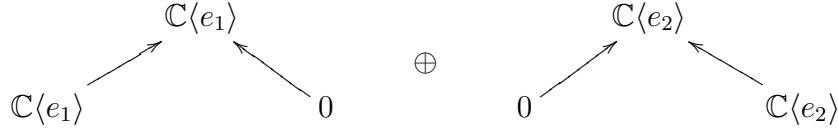
\begin{center}
$\xymatrix @R=0.1pc { &\mathbb C\langle e_1 \rangle\ar@{<-}[ddl]\ar@{<-}[ddr]& & &&\mathbb C\langle e_2\rangle\ar@{<-}[ddl]\ar@{<-}[ddr]&\\
&&& \oplus &&&& \\
\mathbb C \langle e_1\rangle&&0&&0&&\mathbb C \langle e_2 \rangle}$
\end{center}
One can show also that this poset has only finite number of isomorphism
classes of indecomposable objects in $\mathcal P-sp$.
\end{ex}

The \emph{dimension vector} $d=(d_0;d_i)_{i\in \mathcal P}\in \mathbb
N^{|\mathcal P|+1}$ of a given representation $(V;V_i)_{i\in \mathcal
P}$ is defined as $d_0=\dim V$, $d_i=\dim V_i$.
Fixing a dimension vector $d$ we define the variety
$R_{\mathcal P,d}$ of representations of $\mathcal P$ with
dimension $d$
\begin{eqnarray*}
R_{\mathcal P,d}=\left\{(V_i)_{i\in \mathcal P} \in \prod_{i\in \mathcal P} \textrm{Gr}(d_i,d_0)\ \Big |\ V_i\subset V_j, \ i\prec j\right\}.\end{eqnarray*}
The group $GL(d_0)$ acts on $R_{\mathcal P,d}$ via
simultaneous base change $g \cdot  (V_i)_{i\in \mathcal P} = (g(V_i))_{i\in \mathcal P}$, so that the orbits of this action are in one-to-one correspondence with
the isomorphism classes of representations of $\mathcal P$ with
the dimension vector $d$.

\begin{ex} \label{PL4}
    Consider the poset $\mathcal P=(1,1,1,1)$
    consisting of four incomparable elements.
    L.Nazarova \cite{N} and I.Gelfand with V.Ponomarev \cite{GP} completely
    classified indecomposable systems of four subspaces.
    Let us recall the description for the dimension vector $d=(2;1,1,1,1)$:
    \begin{center}
        $
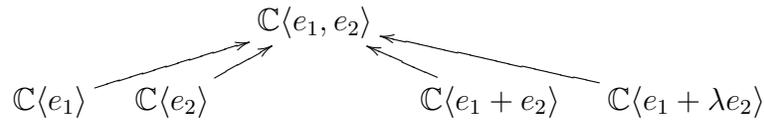
\xymatrix @-1pc { &&\mathbb C\langle e_1,e_2 \rangle&&\\
         \mathbb C \langle e_1\rangle\ar@{->}[urr]&\mathbb C \langle e_2\rangle\ar@{->}[ur]&&\mathbb C \langle e_1+e_2\rangle\ar@{->}[ul]&\mathbb C \langle e_1+\lambda e_2\rangle\ar@{->}[ull]}$
    \end{center}
    It was shown that non-isomorphic indecomposable representations are parametrized by the extended complex plane $\mathbb C \cup \{\infty \}\simeq S^2$
    (when $\lambda=\infty$ the corresponding
    subspace $\mathbb C \langle e_1+\lambda
    e_2\rangle$ is defined to be $\mathbb C \langle e_2\rangle$).
    One can see that when $\lambda=0$, $1$ or
    $\infty$ then the subspaces in corresponding representations are not in general position (this cases are exceptional). Therefore we have the sphere without three points of non-exceptional representations.
\end{ex}

\subsection{Unitary representations. Orthoscalar representations.}

Define the category $\mathcal P-usp$ of unitary representations
of a poset $\mathcal P$.  The objects of $\mathcal P-usp$ are systems
$(U;U_i)_{i\in \mathcal P}$, in which $U$ is a unitary space, each
$U_i$ is a subspace of $U$ and $U_i\subset U_j$ if $i\prec j$. The
set of morphisms between two objects $(U;U_i)_{i\in \mathcal P}$
and $(\tilde U;\tilde U_i)_{i\in \mathcal P}$ consists of
isometric maps $\varphi:U\rightarrow \tilde U$ such that $\varphi(U_i)\subset
U_j$. Two systems are said to be \emph{unitary equivalent} if the
morphism between them can be chosen to be a unitary map. 
An object $(U;U_i)_{i\in \mathcal P}$ is said to be
\emph{indecomposable} it is not unitary equivalent to an orthogonal sum of two
non-zero objects in $\mathcal P-usp$.

\begin{ex} The poset $\mathcal P=(1,1)$ contains a continuous family of unitary non-equivalent representations.
Indeed, consider the representations of the poset $\mathcal P$ as in
Example 2. Taking usual scalar product in $\mathbb C\langle e_1,e_2\rangle $ we view these representations as the objects in $\mathcal P-usp$. For different $\lambda \in (0,1]$ they are indecomposable and unitary non-equivalent.  Moreover, P.Halmos \cite{hal} proved that the representations as in Example 2 gives a complete list of all indecomposable objects in $\mathcal P-usp$. 
\end{ex}

\begin{ex} Consider the poset $\mathcal P=(1,2)$. There exists only a finite number of indecomposable objects in $\mathcal P-sp$ up to the isomorphism. At the same time it is an extremely hard problem  to classify the indecomposable objects  in $\mathcal P-usp$ up to the
unitary equivalence. By S.Kruglyak and Yu.Samoilenko
\cite{KruglyakSamoilenko1,KruglyakSamoilenko2}, such problem
is $*$-wild; that is, it contains the problem of classifying any system of operators in unitary spaces (see \cite[Chapter 3]{OstrovskyiSamoilenko} for more details).
\end{ex}

Fix a weight $\chi=(\chi_0;\chi_i)_{i\in \mathcal P} \in \mathbb
Z^{|\mathcal P|+1}_+$ and consider those objects $(U;U_i)_{i \in \mathcal P}$ in
$\mathcal P-usp$\ that satisfy
\begin{eqnarray*}
    \sum_{i\in \mathcal P} \chi_iP_{U_i}=\chi_0 I,
\end{eqnarray*}
where $P_{U_i}$ is the orthogonal projection onto $U_i$.
Such objects are called \emph{$\chi$-orthoscalar}
(\cite{AlbeverioOstrSamoilenko,KrugNazRoi,KrugRoi}) and form the
subcategory denoted by $\mathcal P_\chi-usp$.
One of the main motivation to study
$\chi$-orthoscalar objects is that any such object gives rise to a polystable 
reflexible sheaves on certain toric variety (see
\cite{Klyachko,KnuS} and references therein).

\begin{ex} Consider the poset $\mathcal P=(1,1)$. The objects in
$\mathcal P_{\chi}$-usp are those systems $(U;U_1,U_2)$ that satisfy
\begin{eqnarray*}
    \chi_1 P_{U_1}+\chi_2 P_{U_2}=\chi_0 I.
\end{eqnarray*}
One can see that the projections $P_{U_1}$ and $P_{U_2}$ commute. Hence the
indecomposable objects are at most one-dimensional. For the
category $\mathcal P_{\chi}$-usp to be non-empty $\chi_0$ should
be equal to $\chi_1+\chi_2$ (trace identity). Moreover, if this
condition is satisfied then the indecomposables in $\mathcal
P_{\chi}-usp$ are in one-to-one correspondence with the
indecomposables in $\mathcal P-sp$.
\end{ex}

\begin{ex} \label{PU4}
Consider the poset $\mathcal P=(1,1,1,1)$. As in Example 3 we
consider the objects with the dimension vector $(2;1,1,1,1)$.
Let us take the weight $\chi=(2;1,1,1,1)$. 
The description of unitary non-equivalent irreducible quadruples of
projections  that satisfy 
\begin{eqnarray*}
    P_{U_1}+P_{U_2}+P_{U_3}+P_{U_4}=2 I,
\end{eqnarray*}
is the following (e.g. \cite[Chapter 2.2]{OstrovskyiSamoilenko}):
\begin{eqnarray*}
    P_{U_1}=\frac{1}{2}\left(
                         \begin{array}{cc}
                           1+a & -b-ic \\
                           -b+ic & 1-a \\
                         \end{array}
                       \right), \quad
    P_{U_2}=\frac{1}{2}\left(
                         \begin{array}{cc}
                           1-a & b-ic \\
                           b+ic & 1+a \\
                         \end{array}
                       \right), \\
    P_{U_3}=\frac{1}{2}\left(
                         \begin{array}{cc}
                           1-a & -b+ic \\
                           -b-ic & 1+a \\
                         \end{array}
                       \right), \quad
    P_{U_4}=\frac{1}{2}\left(
                         \begin{array}{cc}
                           1+a & b+ic \\
                           b-ic & 1-a \\
                         \end{array}
                       \right).
\end{eqnarray*}
Topologically the set of parameters $a, b, c$ is a sphere $S^2$ without
three points where representations split into one-dimensional.
Moreover by \cite{SamMos}, to each triple $(a,b,c)$ one can associate $\lambda \neq 0,1, \infty$
(and vice versa) such that the corresponding system $(U;U_i)_{i\in \mathcal P}$ is
equivalent to a system $(V;V_i)_{i\in \mathcal P}$ which corresponds to $\lambda$ as in Example 3. Hence, the uni-classes of indecomposable
objects in $\mathcal P_{(2;1,1,1,1)}-usp$ with dimension vector
$(2;1,1,1,1)$ are in one-to-one correspondence with the
iso-classes of indecomposable non-exceptional objects in $\mathcal
P-sp$ with dimension vector $(2;1,1,1,1).$
\end{ex}

Recall that $\mathcal P$ is said to be \emph{representation-finite}
(resp. \emph{$*$-representation-finite}) if it has only a finite number of
isomorphism classes of indecomposables in $\mathcal P-sp$
(resp. in $\mathcal P_\chi-usp$ for any weight $\chi$).
P.Gabriel classified all representation-finite quivers together with their
indecomposable representations. Correspondingly, M.Kleiner
\cite{Kleiner1} classified representation-finite posets, and also their
indecomposable representations. There is a direct analogue of the
Gabriel's classification obtained by S.Kruglyak and A.Roiter \cite{KrugRoi}
for unitary orthoscalar representations of quivers (see
\cite{KrugRoi} for the definition of orthoscalarity for quivers). An
analogue of the Kleiner's classification was obtained in \cite{SY}. In particular it turned out that
\begin{eqnarray*}
    \mathcal P\ \ \mbox {is representation-finite}\ \Leftrightarrow     \mathcal P\ \  \mbox {is $*$-representation-finite.}
\end{eqnarray*}

\subsection{Non-commutative Hopf fibration}

Consider the subset of $R_{\mathcal P,d}$ consisting of
the representations which satisfy $\chi$-orthoscalar condition:
\begin{eqnarray*}
    R^{\chi}_{\mathcal P,d}=\left\{ (V_i)_{i \in \mathcal P} \in R_{\mathcal P,d}\ \Big |\ \sum_{i \in \mathcal P} \chi_i P_{V_i}=\chi_0 I
    \right\},
\end{eqnarray*}
where the projections $P_{V_i}$ are taking with respect to the standard
Hermitian metric on $\mathbb C^{d_0}$. The group $U(d_0)$ acts on
$R_{\mathcal P,d}$ as a subgroup of $GL(d_0)$. The orbits of this
action on $R^{\chi}_{\mathcal P,d}$ are in one-to-one
correspondence with the uni-classes of $\chi$-orthoscalar
representations of $\mathcal P$ with dimension vector $d$. 

We will see that the connection between orbit spaces $R_{\mathcal
P,d}/GL(d_0)$ and $R^{\chi}_{\mathcal P,d}/U(d_0)$ is a generalization of 
the following commutative identification.  

The group $\mathbb C^*$
(identified with $GL(1)$) acts on $\mathbb C^n$ by multiplication
$c*z= (c z_1,\ldots,c z_n)$. The corresponding orbit space  is
not Hausdorff (the orbit of $0$ lies in any neighbourhood of any
other orbit), but $\mathbb C^n - \{0\}/\mathbb C^*$ is Hausdorff and homeomorphic to a projective space $\mathbb C \mathbb P^{n-1}$. Consider the subset of $\mathbb C^n$ consisting of those points
$(z_1,\ldots,z_n)$ that satisfy
$z_1\overline{z}_1+\ldots+z_n\overline{z}_n=1$. This subset determines the sphere $S^{2n-1}$. The group $S^1$ (identified with the
unitary group $U(1)$) acts on $S^{2n-1}$ by rotations. The
corresponding orbit space $S^{2n-1}/S^1$ (so called
\emph{Hopf-fibration})
is again a projective space $\mathbb C \mathbb P^{n-1}$. So
we trivially have
\begin{eqnarray*}
\label{Hopf}
    \big\{(z_1,\ldots,z_n)\in \mathbb C^n\ |\ z_1\overline{z}_1+\ldots+z_n\overline{z}_n=1 \big\}/U(1)  \cong  \mathbb C^n - \{0\}/GL(1). 
\end{eqnarray*}

Let $\mathcal P$ be a poset with $n$ non-comparable elements, $d$ be dimension vector and $\chi$ be the weight.  Formally substitute each $z_i$ by the matrix $A_i \in M_{d_0\times d_i}(\mathbb C)$. The corresponding equation between $z_i$ has the following form
\begin{eqnarray*}
    A_1A_1^*+\ldots+A_nA_n^*=I.
\end{eqnarray*}
Take $(\chi^\prime_1,\dots,\chi^\prime_n)=(\chi_1\chi_0^{-1},\ldots,\chi_n\chi_0^{-1})$. Viewing each $\chi_i^\prime$ as $\|A_i\|^2$ we write the non-commutative version of $\overline z_i z_i=|z_i|^2$ as $A_i^*A_i= \chi_i^\prime I_{d_i}$
From last equation we get $A_iA_i^*=\chi_i^\prime P_{{\rm Im}(A_i)}$. 
We say that $(A_1,\ldots,A_n)$ and $(\tilde A_1,\ldots,\tilde A_n)$
lie in the same equivalence class under the action of $U(d_0)$ if there exists $\varphi \in U(d_0)$ such that $\varphi({\rm Im}(A_i))={\rm Im}(\tilde A_i)$ for all $i \in \mathcal P$.  Then the quotient
\begin{eqnarray*}
\left \{(A_i)_{i\in \mathcal P} \in \mathbb(M_{d_0\times d_i}(\mathbb C))_{i\in \mathcal P}\ \ \Big |        
         \begin{array}{c}
                A_1A_1^*+\ldots+A_nA_n^*=I, \\
                A_i^*A_i=\chi_i^\prime I_{d_i}   \\
           \end{array}
           \right\}\Big /U(d_0)
\end{eqnarray*}
parametrizes the equivalence classes of $\chi$-orthoscalar representations with dimension vector $d$ of the poset $\mathcal P$ and can be seen as a \emph{non-commutative Hopf fibration}. On the other hand the quotient 
\begin{eqnarray*}
\left \{(A_i)_{i\in \mathcal P} \in \mathbb(M_{d_0\times d_i}(\mathbb C))_{i\in \mathcal P}\ \ \Big |        
         \begin{array}{c}
                A_i^*A_i=\chi_i^\prime I_{d_i}   \\
           \end{array}
           \right\}\Big /GL(d_0)
\end{eqnarray*}
can be identified with $R_{\mathcal
P,d}/GL(d_0)$. 
In what follows  we show (see Section 4) that for any poset $\mathcal P$ these two quotients are connected similarly to commutative example.

\section{Symplectic quotient.}  \label{s_symplect}

\subsection{Lie groups and algebras. Coadjoint representations}
We briefly recall necessary information about Lie groups, Lie
algebras and coadjoint representations of Lie groups (more
information see, for example, in \cite{KirillovAL}).

By $G$ we denote a Lie group (which is assumed to be
finite-dimensional), $\mathfrak g=\mathfrak{Lie}(G)$ its Lie
algebra (the tangent space to the identity element of $G$) and
$\mathfrak g^*$ denotes the dual Lie algebra. A group $G$ is called
\emph{complex matrix group} if $G$ is a subgroup and at the same
time a smooth submanifold of $GL(n)$. In this paper we mainly
consider the case when $G=U(n)$; that is, the group of unitary
matrices in $n$-dimensional complex space $\mathbb C^n$. The corresponding Lie
algebra $\mathfrak u(n)$ consists of skew-Hermitian matrices.

A Lie group $G$ acts on itself by inner automorphisms
\begin{eqnarray*}A:G\rightarrow Aut(G),\quad A(g):h\mapsto ghg^{-1}.\end{eqnarray*}
Differentiating  we get the \emph{adjoint representation}
$d(A(g))_e=Ad_g:\mathfrak g \rightarrow \mathfrak g$. In the case
when $G$ is a matrix group this representations is given by
\begin{eqnarray*}
    Ad_g(x)=gxg^{-1}, \quad g \in G, \quad x \in \mathfrak{g}.
\end{eqnarray*}
Let $\langle\cdot,\cdot\rangle: \mathfrak g^*\times \mathfrak g
\rightarrow \mathbb K$, $(\xi,x)\rightarrow \langle \xi,
x\rangle=\xi(x)$ be natural pairing between $\mathfrak g$ and
$\mathfrak g^*$. The \emph{coadjoint representation}
$Ad^*:G\rightarrow Aut(\mathfrak g^*)$ is defined by
\begin{eqnarray*}
        \langle Ad^*_g \xi,    x \rangle=\langle \xi,Ad_{g^-1} x \rangle,\quad   x\in \mathfrak g, \quad \xi \in \mathfrak{g^*}.
\end{eqnarray*}
In the case when $G$ is a matrix group the coadjoint
representations is given by
\begin{eqnarray*}
    Ad^*_g(\xi)=g\xi g^{-1}, \quad g \in G, \quad \xi \in \mathfrak{g^*}.
\end{eqnarray*}
Recall that if $\mathfrak g^*$ is semisimple then adjoint and
coadjoint representations are equivalent. By the \emph{coadjoint
orbits} we understand the orbits of $G$ on $\mathfrak{g^*}$.

\begin{ex}
    Let $G=U(n)$. Its dual Lie algebra $\mathfrak{u}(n)^*$
    consists of the
    Hermitian matrices via the identification
    $\mathfrak{u}(n)^*\cong i\mathfrak{u}(n)$.
    Then each coadjoint orbit $\mathcal O_\lambda$ is the set
     Hermitian matrices
    that have the spectrum $\lambda=(\lambda_1,\ldots,\lambda_n)\in \mathbb R^n$.
\end{ex}

We say that a Lie group $G^{\mathbb C}$ is a
\emph{complexification} of $G$, if $G$ is a closed sub-Lie group
of $G^{\mathbb C}$ and $\mathfrak g^{\mathbb
C}=\mathfrak{Lie}(G^{\mathbb C})$ is a vector-space
complexification of $\mathfrak g=\mathfrak {Lie}(G)$; that is,
$\mathfrak g^{\mathbb C}=\mathfrak g \oplus i \mathfrak g$.

\begin{ex} The complexification of $G=SL(n,\mathbb R)$
is the group $SL(n,\mathbb
R)^{\mathbb C}=SL(n,\mathbb C)$, due to obvious identification
$\mathfrak{sl}(n,\mathbb R)^{\mathbb C}=\mathfrak{sl}(n,\mathbb
C)$. The complexification of $G=U(n)$ is $GL(n)$, because $\mathfrak u(n)^{\mathbb C}$ is the full
matrix algebra. 
\end{ex}

\subsection{Symplectic manifolds.}
Let us recall some basic facts about symplectic manifolds (more details see, for example,
\cite{AnaSilva}). Let $M$ be a manifold, and $\omega$ be a closed 2-form on the
tangent space  which assumed to be non-degenerate; that is, for any
tangent vector $x_1 \in T_pM$, there exists another vector $x_2 \in
T_pM$ with $\omega(x_1,x_2)$ nonzero. Due to non-degeneracy of
$\omega$ the manifold $M$ has to be even dimensional.

\begin{ex} Consider the space $\mathbb R^{2n}$ with basis
$\{x_1,\ldots,x_n,y_1,\ldots,y_n\}$ and the form $\omega$ acting by $\omega(x_i,x_j)
=\omega(y_i,y_j)=0$, and $\omega(x_i,y_j) = \omega(y_j,x_i) =
\delta_{ij}$ (Kronecker delta). Such form $\omega$ is called \textit{standard}.
\end{ex}

If $(M,\omega)$ is locally isomorphic to $\mathbb R^{2n}$ with
standard  $\omega$ then $M$ is called a \emph{symplectic
manifold}, and $\omega$ is its \emph{symplectic form}.

\begin{ex} \label{exFubini} Complex projective space $\mathbb C \mathbb P^n$
can be equipped with the Fubini-Study form which makes it into a
symplectic manifold (see, for example, \cite[Section 16]{AnaSilva}). We do not define this form, just note
that one can think of it as the one realized on the quotient space
after the identification $\mathbb C \mathbb P^n=S^{2n+1}/S^1$.
\end{ex}

The coadjoint orbits have a structure of symplectic manifold,
moreover each coadjoint orbit possesses a $G$-invariant symplectic
structure (see \cite[Chapter 1]{Kirillov}). In some cases the opposite is also true: a symplectic Hausdorff
manifold $(M,\omega)$, with the action of a Lie group $G$ which
preserves $\omega$, is isomorphic to a coadjoint orbit of $G$ (see,
for example, \cite{Ko}).

\begin{ex}
    Describe the symplectic form on a coadjoint orbit of $U(n)$.
    It is enough to calculate it in a generic
    point $p \in \mathfrak u(n)^*$.
    Any tangent vector to $p$ looks like
     $ad_p(x)=[p,x]$ for some $x\in \mathfrak u(n)$.
    Then the symplectic form is
    given by $\omega_p(ad_p(x_1),ad_p(x_2))=tr(x_1px_2-x_2px_1)$.
\end{ex}

\begin{ex}
    Consider a coadjoint orbit $\mathcal O_\lambda$ of
    $U(n)$ with $\lambda_1\neq\lambda_2=\lambda_3=\ldots=\lambda_n$.
    The eigenspaces corresponding to  elements
    of $\mathcal O_\lambda$ are
    the line and the orthogonal hyperplane.
    Hence, $\mathcal O_\lambda$ is homeomorphic
    to $\mathbb C \mathbb P^{n-1}$
    and by the Kirillov's construction we have many
     symplectic forms on $\mathbb C \mathbb P^{n-1}$:
     one for each distinct
    pair of real numbers.
    In particular, the Fubiny-Study  form (Example~\ref{exFubini})
    corresponds to the choice $\lambda_1=1$, $\lambda_2=0$.
    The corresponding coadjoint orbit is just the set of
    one-dimensional orthogonal projections.
\end{ex}

\subsection{Moment map}

For a given function $f:M \rightarrow \mathbb C$ its
\emph{symplectic gradient}  $X_f:M\rightarrow T_*M$ (called the
\emph{Hamiltonian vector field}) is defined by the following
equation
\begin{eqnarray*}
    d_{x}f = \omega (x, X_f(p)),
\end{eqnarray*}
where $x$ and $X_f(p)$ are tangent vectors to the point $p\in M$ and
$d_x f$ is the derivation of the function $f$ in the direction $x
\in T_pM$. Since $\omega$ is non-degenerate this defines $X_f$
uniquely.

Let $G$ be a connected Lie group acting on some symplectic manifold
$M$ smoothly preserving the symplectic
form. Then it generates a Lie algebra homomorphism from $\mathfrak g$ to the
Lie algebra of smooth vector fields on $M$. Indeed, an element of
$G$ near the identity gives a diffeomorphism of $M$ close to the
identity; differentiating we get that each tangent vector to the
identity generates a vector field on $M$. By $p \mapsto a_p$ we
denote the vector field on $M$ associated to $a \in \mathfrak g$.

We say that $\mu:M\rightarrow \mathfrak g^*$ is a \emph{moment
map} for the action of $G$ on $M$ if the following two conditions
are satisfied:
\begin{enumerate}
    \item $\mu$ is $G$-equivariant with respect to the action of
    $G$ on $M$ and to the coaction of $G$ on $\mathfrak g^*$; that is, the
    following holds
    \begin{eqnarray*}\mu(g\cdot p)=g\mu(p)g^{-1},\quad p\in M,\ g\in G;\end{eqnarray*}
    \item
    $ d_x \langle \mu(p),a\rangle =\omega(x,a_p),
    \quad \xi \in T_xM$ and $a\in \mathfrak g$. This property means
    that the vector field associated to any $a\in
    \mathfrak g$ equals to the symplectic gradient
    of the function $f_a(p)=\langle
    \mu(p),a\rangle:M\rightarrow \mathbb C$.
\end{enumerate}

If an action has moment map it is said to be \emph{Hamiltonian},
and the corresponding moment map is uniquely determined up to
adding  a constant in $\mathfrak g^*$. On the other hand given a
moment map, one can recover the action of the Lie algebra and the
Lie group on the manifold.

\begin{ex}
    Let $M=\mathbb R^{2}$ and the Lie group $S^1$ acts by
    rotations on $M$.
    It obviously preserves $\omega$, and the corresponding moment map
    $\mu:\mathbb R^{2}\rightarrow \mathbb R$ is defined as
    $\mu(x)=x_1^2+x_{2}^2+{\rm const}$.
\end{ex}

\begin{ex}
    Assume that $M$ is a coadjoint orbit of $G$.
    There exists a unique symplectic structure on $M$
    (called Kirillov-Konstant-Souriau)
    such that the moment map is the embedding
    \begin{eqnarray*}\mu:M \rightarrow \mathfrak{u}(n)^*,\quad \mu:x\mapsto x.\end{eqnarray*}
\end{ex}

\begin{ex} \label{exGrass} Assume the $M=M_{n\times k}(\mathbb C)$
with the symplectic form $\omega(A,B)={\rm tr}(A^*B)-{\rm tr}(B^*A)$, and with
the natural action of $U(k)$ which obviously preserves the
symplectic structure. Then  $\mu:M\rightarrow \mathfrak u(k)^*$
 is given by $\mu:A\mapsto \frac{i}{2}AA^*+{\rm const}\cdot I$.
\end{ex}

An important property of symplectic manifolds is that having two
manifolds $(M_1,\omega_1)$ and $(M_2,\omega_2)$ one can form the
product manifold $M_1\times M_2$ with symplectic structure $\pi^*_1
\omega_1+\pi^*_2 \omega_2$, where $\pi_i:M_1\times M_2 \rightarrow
M_i$ is the projection onto $i$-th factor, and $\pi_i^*$ is the
pull-back of $\pi_i$.

Assume that $M_1$ and $M_2$ are symplectic manifolds with the action
of the same Lie group $G$ and corresponding moment maps $\mu_1$ and
$\mu_2$. Then the symplectic manifold $M_1\times M_2$ 
possesses the diagonal action of $G$ with the moment map
\begin{eqnarray*}\mu:M_1\times M_2 \rightarrow \mathfrak{g}^*,\qquad \mu:(x,y)\mapsto \mu_1(x)+\mu_2(y).\end{eqnarray*}

\begin{ex} \label{exHopf}
    Let $M=\mathbb R^{2n}$ with the action of $S^1$. Then
    $\mu:\mathbb R^{2n}\rightarrow \mathbb R$ is defined as
    $\mu(x)=x_1^2+\ldots+x_{2n}^2+{\rm const}$.
\end{ex}

\begin{ex} \label{exPos}
    Consider the product $\mathcal O_{\lambda^{(1)}} \times \ldots \times \mathcal O_{\lambda^{(m)}}$
    of coadjoint orbits of $U(n)$. Moment map takes the set
    of Hermitian matrices $(A_1,\ldots,A_m)$
    with $\sigma(A_i)=\lambda^{(i)}$ to their sum:
    \begin{eqnarray*}
        \mu:(A_1,\ldots,A_m) \mapsto A_1+\ldots+A_m+{\rm const}\cdot I.
    \end{eqnarray*}
\end{ex}

\subsection{Symplectic quotient}

Assume that $M$ is a manifold with Hamiltonian action of a Lie group
$G$ and with the corresponding moment map $\mu:M\rightarrow g^*$.
Then $G$ acts on the fiber $\mu^{-1}(0)\subset M$.  The following
theorem says that the corresponding orbit space has the structure of
a symplectic manifold.

\begin{thm} (Marsden-Weinstein)
The quotient space $\mu^{-1}(0)/G$ is a symplectic manifold. If the
action of $G$ is free on $\mu^{-1}(0)$ then $\mu^{-1}(0)/G$ has
dimension $\dim M - 2\dim G$.
\end{thm}

The manifold $\mu^{-1}(0)/G$ is called the \emph{symplectic
quotient}.

\begin{ex}
    Consider the manifold and moment map
    as in Example~\ref{exHopf}.
    The constant in the moment map can be chosen to be $-1$,
    hence $\mu^{-1}(0)$
    is the sphere $S^{2n-1}$ and the symplectic quotient is
    homeomorphic to
    $\mathbb C \mathbb P^{n-1}$ (Hopf fibration).
\end{ex}

\begin{ex}
Let $M$ be the manifold as in Example~\ref{exGrass}. We choose the
constant in the moment map to be equal to $-\frac{i}{2}$. Then
$\mu^{-1}(0)=\{ A \in M_{n\times k}(\mathbb C)\ |\ AA^*=I\}$ and the
quotient $\mu^{-1}(0)/U(k)$ is the set of
$k$-dimensional subspaces in $n$-dimensional space; that is,
$\textrm{Gr}(k,n)$.
\end{ex}

\section{Quick GIT. Stability conditions and Kirwan-Ness theorem} \label{s_git}

Assume that $M$ is a variety and $G$ is a linear
algebraic group acting on $M$. Group $G$ is assumed to be
\emph{reductive}; that is, it is a complexification of some compact Lie group.
Geometric Invariant Theory  tries to build the quotient of $M$ by
$G$. The main problem is that the quotient $M/G$ may not exist in
the category of algebraic or projective varieties (especially when the group $G$
is not finite). One of the possible solution is to remove
some points from $M$, by taking an open subset $M^\prime$ of $M$ as large as
possible such that $M^\prime/G$ is a variety. The following elementary example explains the idea. The space $\mathbb C\mathbb P^n$ can be described as the GIT quotient of the space
$\mathbb C^{n+1}$ by $\mathbb C^*$, namely $\mathbb C \mathbb
P^{n}=(\mathbb C^{n+1}-0)/ \mathbb C^*$.

\subsection{Projective spectrum and GIT quotients}

Let $R=\bigoplus_{k \in \mathbb N} R_k$ be a \emph{graded} ring
(the product of an element from $R_n$ with an element from $R_m$
lies in $R_{n+m}$). The \emph{projective spectrum} of the ring $R$
(denoted by $\Pj(R)$) is defined as the set of maximal graded
ideals of $R$. There is alternative definition in terms of
ordinary spectrum of the ring. We have the projection
$R\twoheadrightarrow R_0$, hence $\Sp(R_0)\hookrightarrow \Sp(R)$
($\Sp(\cdot)$ is contravariant). The group $\mathbb C^*$ acts on
$R$ by rotating $R_k$'s. Then
\begin{eqnarray*}
\textrm{Proj}(R)=(\textrm{Spec}(R) \setminus \textrm{Spec}(R_0))/
\mathbb C^*.
\end{eqnarray*}

\begin{ex}
    Let $R=\mathbb
    C[x_1,\ldots,x_n]$. Then $R_0=\mathbb C$. Hence,
    $\Pj(\mathbb C[x_1,\ldots,x_n])=(\Cn^n\setminus \{0\})/\mathbb
    C^*$ which is identified with $\mathbb C\mathbb P^{n-1}$.
\end{ex}

A \emph{polarization} $L$ of $M$ is a line bundle $L\rightarrow M$, equipped with 
the action of $G$. Graded ring $R(M)$ associated to $M$ is defined as
\begin{eqnarray*}
    R(M)=\bigoplus_{n\geq 0} \Gamma(M,L^n),
\end{eqnarray*}
where $L^n$ is the $n$-th tensor power of $L$ and $\Gamma(M,L^n)$ is the set of $G$-invariant sections of $L^n$. The \emph{GIT quotient} of $M$ by $G$ with respect to the line
bundle $L$ is defined as $\Pj(R(M))$  and denoted by $M//G$.

\medskip 
A point $p\in M$ is:
\begin{itemize}
\item \emph{semistable} if $s(p)>0$ for some section $s\in
\Gamma(M,L^n)$, in which $n>0$; \item \emph{polystable} if $p$ is
semistable and the orbit $\{ g \cdot  p\ |\ g \in G\}$ is closed; \item \emph{stable} if
$p$ is polystable with a finite stabilizer;
\item \emph{unstable}
if $m$ is not semistable.
\end{itemize} Let $M^{ss}$,
$M^{ps}$, $M^s$ and $M^{us}$ be  semistable, polystable, stable
and unstable locus respectively. D.Mumford gave a numerical
criterion to decide the type of the point $m \in M$ and proved 
that the GIT quotient $M//G$ can be identified with $M^{ps}/G$
(see \cite[Chapter 2]{MumfordForganyKirwan} for the details).

\begin{ex}(\cite[Chapter 8]{MumfordForganyKirwan}) \label{exPoints}
 Let  $M=(\mathbb C \mathbb P^1)^n$ with $G=SL(2)$ acting diagonally.
 Then we have:
\begin{eqnarray*}        
		M^{s}&=&\{(m_1,\ldots,m_n)\in (\mathbb C \mathbb P^1)^n\ |\ \mbox{at most n/2
        points equal}\},\\
        M^{ss}&=&\{(m_1,\ldots,m_n)\in (\mathbb C \mathbb P^1)^n\ |\ \mbox{less then n/2
        points equal}\},\\
        M^{ps}-M^{s}&=&\{(m_1,\ldots,m_n)\in (\mathbb C \mathbb P^1)^n\ | \  \#\{m_1,\ldots,m_n\}=2\}.
\end{eqnarray*}
\end{ex}

\subsection{Kirwan-Ness Theorem}
Assume that $M$ is a symplectic manifold with an action of a compact group $G$, which 
preserves symplectic form and has the moment map $\mu$. Then there
exists a unique continuation of the action of $G$ to the action of
$G^\mathbb C$ (its complexification). By \cite{MumfordForganyKirwan}, there exists an inclusion of $\mu^{-1}(0)$ into $M^{ps}$. The following theorem is a
fundamental fact (proved by Kirwan and Ness independently)
connecting the symplectic  and the GIT quotients.

\begin{thm} (Kirwan-Ness)
The inclusion $\mu^{-1}(0)$ into $M^{ps}$ induces a homeomorphism
\begin{eqnarray*}
    \mu^{-1}(0)/G \cong M//G^{\mathbb C}.
\end{eqnarray*}
\end{thm}

\begin{ex}
    Let $M=\mathbb R^{2n}$, and let $S^1$ acts by rotations on $M$.
    The complexification of $S^1$ is $\mathbb C^*$
    which acts by multiplication on $M$.
    The Kirwan-Ness theorem is just an identification
    between $S^{2n-1}/S^1$ and $\mathbb C \mathbb P^{n-1}$.
\end{ex}

\begin{ex}
    Consider the space $M=(\mathbb C \mathbb P^1)^n$ with diagonal action of $G=SL(2)$. By  Example~\ref{exPoints}, $M^{s}=M^{ss}$ if and only if $n$ is odd. Take the maximal
compact subgroup $SU(2)$ of $G$. Its dual Lie algebra can be
identified with the real vector space $\mathbb R^3$. The
corresponding moment map takes the $n$-tuple of points on unit
sphere to their sum in $\mathbb R^3$. Then the fiber $\mu^{-1}(0)$
consists of the points with the center of gravity in the origin.
It is not hard to check directly that $\mu^{-1}(0)\subset M^{ss}$.
In the case when $n=4$ it can be straightforwardly proved that 
$(\mathbb C \mathbb P^1)^n//SL(2)$ and $\mu^{-1}(0)/SU(2)$ are isomorphic to 
a sphere $\mathbb C \mathbb P^1$ without three removed points. We refer the reader to
\cite[Chapter 8]{MumfordForganyKirwan} for the details.
\end{ex}

Note that the correspondence in previous example is the same as the corres\-pon\-dence
between Example~\ref{PL4} and Example~\ref{PU4}.

\section{Correspondence between unitary and linear representations
of the posets} \label{s_corr}

For a given poset $\mathcal P$ fix a dimension vector $d=(d_0;d_i)_{i\in \mathcal P}$ and a weight $\chi=(\chi_0;\chi_i)_{i\in \mathcal P}\in \mathbb N^{|\mathcal P|+1}$, such that $\sum_{i\in \mathcal P}\chi_i d_i=\chi_0 d_0$.  Let $\lambda^{(i)}=(\underbrace{\chi_i,\dots,\chi_i}_{d_i},\underbrace{0,\dots,0}_{d_0-d_i})$, $i\in\mathcal P$.  The variety $R_{\mathcal P,d}$ is a subset of the product
of Grassmanians $\prod_{i\in \mathcal P} \textrm{Gr}(d_i,d_0)$. Given an Hermitian metric on $\mathbb C^{d_0}$ one can view each $\textrm{Gr}(d_i,d_0)$ as the set of $d_0\times d_0$ matrices $\chi_i P_i$,
where $P_i$ is the $d_i$-dimensional orthoprojection. Therefore $\textrm{Gr}(d_i,d_0)$ is identified with the coadjoint orbit $\mathcal O_{\lambda^{(i)}}$. Our first aim is to show how to embed the variety $R_{\mathcal P,d}$ into some projective space so that its 
symplectic structure coincides with symplectic structure on the product of corresponding coadjoint orbits. We use slightly modified standard construction (see, for example, \cite[Chapter 11]{Dolgachev} and \cite{Knu}).

    A standard way  to embed $\textrm{Gr}(d_i,d_0)$ into a projective space
    is via Plucker embedding; that is, for an element
    $V_i \in \textrm{Gr}(d_i,d_0)$ we take its basis vectors $a_j$ and
    wedge them together $a_1 \wedge\cdots\wedge a_{d_i}$
    obtaining an element of
    $\mathbb P(\wedge^{d_1}\mathbb C^{d_0})$
    (this is well-defined because if we change a basis then
    $a_1 \wedge\cdots\wedge a_{d_1}$ changes by a scalar).

    If $V$ is some vector space then one can form a
    symmetric tensor $d$-power of $V$ denoted
    $Sym^d(V)$ (for $d<0$ one takes $Sym^{-d}(V^*)$, for
    $d=0$ we have  $Sym^0(V)=\mathbb C$).
    We take the symmetric tensor $\chi_i$-power of the
    space $\wedge^{d_i}\mathbb
    C^{d_0}$. By Veronese map, we embed the projective space $\mathbb P (V)$
    into the space $\mathbb P({Sym^d(V)})$. Hence, we have the
    following sequences of inclusions:
    \begin{eqnarray*}
        \textrm{Gr}(d_i,d_0) \hookrightarrow \mathbb P(\wedge^{d_i}\mathbb
        C^{d_0}) \hookrightarrow \mathbb P ({Sym^{\chi_i}(\wedge^{d_i}\mathbb
        C^{d_0})}).
   \end{eqnarray*}
    Now it is a routine to check that the symplectic form on
    $\textrm{Gr}(d_i,d_0)$ (taken as the restriction of the
     Fubiny-Study form on
     $\mathbb P({Sym^{\chi_i}(\wedge^{d_i}\mathbb
        C^{d_0})})$)  coincides with the symplectic form on the
        corresponding coadjoint orbit $\mathcal O_{\lambda^{(i)}}$.

  Correspondingly for the product of Grassmanians
  $\prod_{i \in \mathcal P} \textrm{Gr}(d_i,d_0)$ we have the embedding
  \begin{eqnarray*}
        \prod_{i \in \mathcal P} \textrm{Gr}(d_i,d_0)
        \hookrightarrow \prod_{i \in \mathcal P}
        \mathbb P(Sym^{\chi_i}(\wedge^{d_i}\mathbb C^{d_0})).
  \end{eqnarray*}
  Using the Segre map $\mathbb P^n\times \mathbb P^m \hookrightarrow \mathbb
  P^{(n+1)(m+1)-1}$ we embed the last product into
  \begin{eqnarray*}
     \mathbb P\left( \prod_{i \in \mathcal P} Sym^{\chi_i}(\wedge^{d_i}\mathbb C^{d_0})\right).
 \end{eqnarray*}
 
Consider the moment map of the action of
$U(d_0)$ on $\prod_{i \in \mathcal P} \textrm{Gr}(d_i,d_0)$ after the embed\-ding. The Fubiny-Study form on $\mathbb P\left( \prod_{i
\in \mathcal P} Sym^{\chi_i}(\wedge^{d_i}\mathbb C^{d_0})\right)$ coincides (as we
mentioned above) with the symplectic form on the product of coadjoint
orbits $\prod_{i \in \mathcal P} \mathcal O_{\lambda^{(i)}}$. Knowing the form of the moment map 
$\mu:(P_i)_{i \in \mathcal P} \mapsto \mathfrak u(d_0)^*$ 
(see Example~\ref{exPos}) and taking a constant in the moment map to
be equal to $-\chi_0$ we have the following
\begin{eqnarray*}
  \mu^{-1}(0)= \left \{ (P_i)_{i \in \mathcal P}\in (M_{d_0}(\mathbb C))_{i \in \mathcal P}\,  \scalebox{1.5}{\Big |}\!
    \begin{array}{c}
        P_i=P_i^*=P_i^2,\ \textrm{rank}(P_i)=d_i, \\
        P_iP_j=P_iP_j=P_i, \quad i\prec j, \\
        \sum_{i\in \mathcal P}\chi_i P_i=\chi_0I   
     \end{array}  \right\}.
\end{eqnarray*}

The continuation of the action of $U(d_0)$ to its complexification $GL(d_0)$ coincides with the action of $GL(d_0)$ on $R_{\mathcal P,d}$. Using Muformd's numerical criterion it was shown in \cite[Chapter 11]{Dolgachev} and \cite[Theorem 2.2]{hu} the set of stable representations
with respect to induced action of $GL(d_0)$ consists of those representations
$(V;V_i)_{i\in\mathcal P}$ that satisfy
\begin{eqnarray*}
    \frac{1}{\dim K} \sum_{i \in \mathcal P}\chi_i \dim(V_i \cap K)
    <\frac{1}{\dim V} \sum_{i \in \mathcal P}\chi_i \dim V_i
\end{eqnarray*}
for each proper subspace $K \subset V$.

Straightforward calculations  show that if $(P_i)_{i\in \mathcal P} \in \mu^{-1}(0)$ then the system of subspaces $(\mathbb C^{d_0};{\rm Im} P_i)_{i\in \mathcal P}$ is \textit{$\chi$-polystable}; that is, it decomposes into a direct sum of stable representations $(\tilde V;\tilde V_i)_{i\in \mathcal P}$ which satisfy $\sum_{i
\in \mathcal P}\chi_i\dim \tilde V_i=\chi_0 \dim \tilde V_0$. Denote by $\Phi$ the corresponding map from $Ob(\mathcal P_\chi-usp)$ to $Ob(\mathcal P-sp)$. As a consequence of Kirwan-Ness theorem we have the following:

\begin{thm} \label{posCor}
Let $\mathcal P$ be a poset,  $(d_0;d_i)_{i \in \mathcal P}$ be a dimension vector and $\chi=(\chi_0;\chi_i)_{i\in \mathcal P}$ be a weight such that $\sum_{i\in \mathcal P} \chi_i d_i=\chi_0 d_0$. Then $\Phi$ induces 
a homeomorphism (with respect to the usual orbit-space
topology) between
\begin{eqnarray*}
   \left \{ (P_i)_{i \in \mathcal P} \in (M_{d_0}(\mathbb C))_{i \in \mathcal P}\ \scalebox{1.5}{\Big |}
    \begin{array}{c}
        P_i=P_i^*=P_i^2,\ \textrm{rank}(P_i)=d_i,  \\
        P_iP_j=P_iP_j=P_i, \quad i\prec j, \\ \quad
        \sum_{i\in \mathcal P}\chi_i P_i=\chi_0I
    \end{array}  \right\} \Big / U(d_0)
\end{eqnarray*}
and
\begin{eqnarray*}
   \left \{(V;V_i)_{i \in \mathcal P}\in R_{\mathcal P,d} \ \ \Big |\
    \begin{array}{c}
        (V;V_i)_{i \in \mathcal P}\ is\ \chi-polystable
    \end{array}
     \right \} \Big / GL(d_0).
\end{eqnarray*}
\end{thm}
Let us show some consequences of this theorem.
\begin{enumerate}
    \item In \cite{KrugNazRoi} it was shown that if two objects in
$\mathcal P_{\chi}-usp$ are equivalent as the objects in $\mathcal
P-sp$  then they are equivalent in $\mathcal P_{\chi}-usp$. This follows as a consequence of Theorem~\ref{posCor}.   
 \item In \cite{KrugNazRoi} it was shown
    that each indecomposable object in
$\mathcal P_{\chi}-usp$ is Schurian; that is, it has trivial endomorphism ring. In fact this result can be deduced from the previous observations. Indeed, each indecomposable object in $\mathcal
P_{\chi}-usp$ corresponds to some $\chi$-stable object in
$\mathcal P-sp$. But stable objects are Schurian.

\end{enumerate}

\section{Calculating the dimensions of the quotients} \label{s_dimension}

\subsection{Bound quivers}

For a finite  quiver $Q$ denote by $Q_0$ the set of its vertices and by $Q_1$ the set of
its arrows denoted by $\rho:i\rightarrow j$ for $i,j\in Q_0$. In
the following we only consider quivers without oriented cycles. A
finite-dimensional \textit{representation} of $Q$ is given by a
tuple
\begin{eqnarray*}X=((X_i)_{i\in Q_0},(X_{\rho})_{\rho \in Q_1}:X_i\rightarrow X_j)\end{eqnarray*}
of finite-dimensional vector spaces and linear maps between them. We say that $X$ is \emph{subspace representation} if all
maps $X_{\rho}$ are injective. The dimension vector
$\underline{\dim}X\in\mathbb{N}^{|Q_0|}$ of $X$ is defined by
\begin{eqnarray*}\underline{\dim}X=(\dim X_i)_{i\in Q_0}.\end{eqnarray*}
The variety $R_{Q,d}$ of representations of $Q$ with
the dimension vector $d \in \mathbb{N}^{|Q_0|}$ is defined as the
affine complex space
\begin{eqnarray*}R_{Q,d}=\bigoplus_{\rho:i\rightarrow j} \mathrm{Hom}( \mathbb C^{d_i},\mathbb C^{d_j}).\end{eqnarray*}
The algebraic group $GL_d=\prod_{i\in Q_0} GL(d_i)$ acts on
$R_{Q,d}$ via simultaneous base change; that is,
\begin{eqnarray*}(g_i)_{i\in Q_0}\cdot 
(X_{\rho})_{\rho\in Q_1}=(g_jX_{\rho}g_i^{-1})_{\rho:i\rightarrow
j}.\end{eqnarray*} The orbits of this action are in a bijection with the isomorphism classes
of representations of $Q$ with the dimension vector
$d$.

Let $\mathbb C Q$ be the path algebra of $Q$ (see, for example, \cite[Chapter 2]{ass} for the definition)  and let $RQ$ be the
arrow ideal. A relation in $Q$ is a $\mathbb C$-linear combination
of paths of length at least two which have the same head and tail.
For a set of relations $(r_j)_{j\in J}$ we can consider the
admissible ideal $I$ generated by these relations, that means that
$RQ^m\subseteq I\subseteq RQ^2$ for some $m\geq 2$. We say
that a representation $X$ of $Q$ is \emph{bound} by $I$ (or
a representation of the \emph{bound quiver} $(Q,I)$) if
$X_{r_j}=0$ for all $j\in J$. For every dimension vector this
defines a closed subvariety of $R_{Q,d}$ denoted by $R_{(Q,I),d}$.
Fixing a minimal set of relations generating $I$, we
denote by $r(i,j,I)$ the number of relations with starting vertex
$i$ and terminating vertex $j$. Following \cite{bon}, for the
dimension of $R_{(Q,I),d}$ we have the following lower bound
\begin{eqnarray*}\dim R_{(Q,I),d}\geq\dim R_{Q,d}-\sum_{(i,j)\in (Q_0)^2}r(i,j,I)d_i d_j.\end{eqnarray*}

Let $C_{(Q,I)}$ be the Cartan matrix of $(Q,I)$; that is, $c_{j,i}=\dim
e_i(\mathbb C Q/I)e_j$ where $e_i$ denotes the primitive idempotent
(resp. the trivial path) corresponding to the vertex $i$. On $\Zn^
{|Q_0|}$ a non-symmetric bilinear form  is
defined by
\begin{eqnarray*}\Sc{d}{{\tilde d}}:=d^t(C_{(Q,I)}^{-1})^t{\tilde d}.\end{eqnarray*}
For two representation $X$ and $Y$ we have
the following homological interpretation of this form:
\begin{eqnarray*}\Sc{X}{Y}:=\Sc{\underline{\dim} X}{\underline{\dim} Y}=\sum_{i=0}^{\infty}(-1)^i\dim\Ext^i(X,Y).\end{eqnarray*}
If $Q$ is unbound, for two representations $X$, $Y$ of $Q$ with
$\underline{\dim} X=d$ and $\underline{\dim} Y={\tilde d}$ we have $\Ext^i(X,Y)=0$ for $i\geq 2$ and 
\begin{eqnarray*}\Sc{X}{Y}=\dim\Hom(X,Y)-\dim\Ext(X,Y)=\sum_{q\in Q_0}d_q {\tilde d}_q-\sum_{\rho:i\rightarrow j\in Q_1}d_i{\tilde d}_j.\end{eqnarray*}

\subsection{Bound quiver and posets}

Recall briefly the relation between
representations of posets and representations of bound quivers.

Let $\mathcal P$ be a poset and $Q^{\mathcal P}$ its Hasse quiver.
All arrows of $Q^{\mathcal P}$  are oriented to one vertex $q_0$
which is called the root. Let $d \in\Nn^{|\mathcal P|+1}$ be a
dimension vector. By $S_{Q^{\mathcal P},d}\subset R_{Q^{\mathcal P},d}$ we denote the open subvariety of subspace representations. For
every (non-oriented) cycle
$\rho_1\ldots\rho_n\tau_k^{-1}\ldots\tau_1^{-1}$ with
$\rho_i,\,\tau_j\in Q_1$ and $\rho_i\neq\tau_j$ we define a relation
\begin{eqnarray*}r=\rho_1\ldots\rho_n-\tau_1\ldots\tau_k.\end{eqnarray*} Let $I$ be the
ideal generated by all such relations.

Let $\pi=(V;V_i)_{i\in \mathcal P}$ be a representation of
$\mathcal P$ with the dimension vector $d$. This defines a
representation $X(\pi)\in S_{(Q^{\mathcal P},I),d}$ satisfying the
stated relations. Indeed, every inclusion $V_i\subset V_j$ defines
an injective map $X(\pi)_{\rho:{i\rightarrow j}}:V_i\rightarrow V_j$. Thus it
defines an element of $S_{(Q^{\mathcal P},I),d}$. For two
arbitrary representations $\pi=(V;V_i)_{i\in \mathcal P}$ and $\pi'=(W;W_i)_{i\in \mathcal P}$ a
morphism $f:\pi\rightarrow \pi'$ induces a morphism
$X(f):X(\pi)\rightarrow X(\pi')$, where
$X(f)_i:=f|_{V_i}:X(V)_i\rightarrow X(W)_i$.

Conversely, let $X\in S_{(Q^{\mathcal P},I),d}$. This gives rise
to a representation $\pi(X)$ of $\mathcal{P}$ by defining
$\pi(X)_q=X_{\rho^q_n}\circ\ldots\circ X_{\rho_1^q}(X_q)$ for some
path $p_q=\rho_1^q\ldots\rho_n^q$ from $q$ to $q_0$. This
definition is independent of the chosen path. Moreover, every
morphism $\varphi=(\varphi_q)_{q\in Q_0}:X\rightarrow Y$ defines a
morphism $\pi(\varphi)$ which is induced by
$\varphi_{q_0}:X_0\rightarrow Y_0$.

Therefore we get an equivalence between the categories of representations
of $\mathcal{P}$ and subspace representations of $Q^{\mathcal P}$
bound by $I$. This equivalence also preserves dimension vectors.

\subsection{GIT quotient of quivers}
A representation $X$ of $(Q,I)$ is \emph{semistable} (resp.
\emph{stable}) with respect to a form  $\Theta\in
\mathrm{Hom}_{\mathbb{Z}}(\mathbb{Z}^{|Q_0|},\mathbb{Z})$ if  $\Theta(\underline{\dim}X)=0$ and 
for all subrepresentations  $Y\subset X$ (resp.
all proper subrepresentations $0\neq Y\subsetneq X$) we have:
\begin{eqnarray*} 
  \Theta(\underline{\dim}Y) \geq 0 \quad  \mbox{ (resp. } \Theta(\underline {\dim} Y)>0 ).
\end{eqnarray*}
Denote the set of semistable (resp. stable) points by
$R^{ss}_{(Q,I),d}$ (resp. $R^s_{(Q,I),d}$). A.King in \cite{kin} built 
the GIT quotient $R_{(Q,I),d}//GL_{d}$ for quivers. This quotient parametrizes
polystable representations. For a stable representation $X$ we have that its orbit is of maximal
possible dimension. The scalar matrices act trivially on
$R_{(Q,I),d}$, hence the isotropy group is one-dimensional. Therefore, if the
quotient $R_{(Q,I),d}//GL_{d}$ is not empty, for its dimension we
have the lower bound
\begin{eqnarray*}\dim R_{(Q,I),d}//GL_{d}&=&\dim R_{(Q,I),d}-(\dim G_d-1)\\&\geq& 1-\sum_{i\in Q}d_i^2+\sum_{\rho:i\rightarrow j\in Q_1}d_id_j-\sum_{(i,j)\in (Q_0)^2}r(i,j,I)d_id_j.
\end{eqnarray*}
In the case $I=0$ we have
$\dim R_{(Q,I),d}//GL_{d}=1-\langle d, d\rangle$.

\subsection{Dimension of symplectic quotients of the posets}

Assume that the representation $\pi=(V;V_i)_{i\in \mathcal P}$ is
stable with the weight $(\chi_0;\chi_i)_{i\in \mathcal P} \in \mathbb Z^{|\mathcal P|+1}$; that is 
\begin{eqnarray*}
    \frac{1}{\dim K} \sum_{i \in \mathcal P}\chi_i \dim(V_i \cap K)
    <\frac{1}{\dim V} \sum_{i \in \mathcal P}\chi_i \dim V_i
\end{eqnarray*}
for each proper subspace $K \subset V$ and 
$    \sum_{i \in \mathcal P}\chi_i \dim V_i=\chi_0\dim V$.
Then the corresponding representation $X(\pi)$ of the bound quiver $Q^{\mathcal P}$ is
stable with respect to the form $\Theta=(\chi_0 ;-\chi_i )_{i\in \mathcal P}$
and vise versa. Therefore, using Theorem~\ref{posCor} we have that if there exists at least one
object in $\mathcal P_\chi-usp$ with dimension vector $(d_0;d_i)_{i\in \mathcal P}$
then the unitary classes of all objects with this dimension vector
depend on at least
\begin{eqnarray*}
1-\sum_{i\in Q^{\mathcal P}_0}d_i^2+\sum_{\rho:i\rightarrow
j\in Q^{\mathcal P}_1}d_id_j-\sum_{(i,j)\in (Q^{{\mathcal P}}_0)^2} r(i,j,I)d_id_j.
\end{eqnarray*}
complex parameters. When the poset is primitive we have that the
number of parameters is $ 1-\sum_{i\in Q^{\mathcal
P}_0}d_i^2+\sum_{\rho:i\rightarrow j\in Q^{\mathcal
P}_1}d_id_j=1-\langle d,d \rangle$.

\begin{ex} Consider the case when $\mathcal P$ is primitive and representation-tame. Hence the corresponding quiver $Q^{\mathcal P}$ is
one of the affine Dynkin quivers; that is, $Q^{\mathcal P}=\tilde
D_4$, $\tilde E_6$, $\tilde E_7$ or $\tilde E_8$. Let us take the
dimension vector $d_{\mathcal P}$ to be the minimal imaginary root
of  $Q^{\mathcal P}$:
\begin{eqnarray*}
    d_{(1,1,1,1)}&=&(2;1,1,1,1),\\
    d_{(2,2,2)}&=&(3;1,2;1,2;1,2),\\
    d_{(1,3,3)}&=&(4;2;1,2,3;1,2,3),\\
    d_{(1,2,5)}&=&(6;3;2,4;1,2,3,4,5).
\end{eqnarray*}
Using Theorem~\ref{posCor} and the dimension formula we have at most $1-\langle d_{\mathcal
P}, d_{\mathcal P} \rangle=1$--parametric family of uni-classes of
indecomposable objects in $\mathcal P_\chi-usp$ with dimension $d_{\mathcal
P}$. One can check
that it is always possible to choose the weight $\chi_\mathcal{P}$
such that there exists at least one stable representation, hence for these
weights the set of uni-classes of orthoscalar representations of $\mathcal P$ in dimension
$d_{\mathcal P}$ is isomorphic to $\mathbb C\mathbb P^1$ (because the GIT quotient is one-dimensional and rational).
\end{ex}

\begin{ex} Consider representation-tame non-primitive critical poset
$\mathcal P=(N,4)$, which has the following Hasse quiver
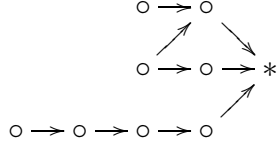
\begin{center}
$\xymatrix @-1pc
            {&&{\circ}\ar@{->}[r]&{\circ}\ar@{->}[dr]\\
             &&{\circ}\ar@{->}[r]\ar@{->}[ur]&{\circ}\ar@{->}[r]&{*}\\
                {\circ}\ar@{->}[r]&{\circ}\ar@{->}[r]&{\circ}\ar@{->}[r]&{\circ}\ar@{->}[ur]}$
\end{center}
Fix the dimension vector $d_{(N,4)}=(5;2,4,3,2;1,2,3,4)$. One
can check that with respect to the weight
$\chi=(5;2,1,1,2;1,1,1,1)$ there exists at least one stable
representation in dimension $d_{(N,4)}$, hence the corresponding GIT quotient is
not empty. Due to the dimension formula we have at least
$1-\langle d_{\mathcal P},d_{\mathcal P} \rangle+2\cdot5=1$--parametric family of uni-classes of representation of
$\mathcal P_\chi-usp$. Moreover, it is possible to show that this
is a maximal number of parameters. Therefore, we have that for
the weight $\chi=(5;2,1,1,2;1,1,1,1)$ the set of uni-classes of
orthoscalar representations of $\mathcal P$ in dimension
$d_{(N,4)}$ is isomorphic to $\mathbb C\mathbb P^1$.
\end{ex}

\ack
This work was done during the visit of the third author to the
University of S\~ao Paulo as a postdoctoral fellow. This author is grateful to the University of S\~ao Paulo for hospitality  and to  Fapesp
for financial support (2010/15781-0).
 The first author is supported in part by  CNPq 
(301743/2007-0) and by  Fapesp  (2010/50347-9).

\end{document}